\def\A{\leavevmode\setbox0\hbox{A}\lower1.4ex\hbox 
to\wd0{\hss`}\kern-.9\wd0A}
\def\E{\leavevmode\setbox0\hbox{E}\lower1.4ex\hbox 
to\wd0{\hss`\/}\kern-.9\wd0E}
\def\a{\leavevmode\setbox0\hbox{a}\lower1.4ex\hbox 
to\wd0{\hss`\/}\kern-\wd0a}
\def\e{\leavevmode\setbox0\hbox{e}\lower1.4ex\hbox 
to\wd0{\hss`\/}\kern-\wd0e}
\newcommand{\be}{\begin{equation}}
\newcommand{\ee}{\end{equation}}
\newcommand{\ba}{\begin{array}}
\newcommand{\ea}{\end{array}}
\newcommand{\beqn}{\begin{eqnarray}}
\newcommand{\eeqn}{\end{eqnarray}}
\newcommand{\nn}{\nonumber}
\newcommand{\p}{\prime}

\font\symb=msam7
\def\znakr{\raise1.5pt\hbox{\symb\char66\kern-2pt\char74}}
\def\znakl{\raise1.5pt\hbox{\symb\char73\kern-2pt\char67}}

\def\normalsize{
\setlength{\textheight}{23cm}
\setlength{\textwidth}{15cm}
\setlength{\topmargin}{-2.0cm}
\setlength{\hoffset}{-0.5cm}
\setlength{\leftmargin}{-1cm}
\setlength{\rightmargin}{2.0cm}}

\documentstyle[12pt]{article}
\normalsize

\begin{document}
\title{Two-Dimensional Centrally Extended Quantum Galilei Groups
and their Algebras}
\author{ Anna Opanowicz\thanks{Supported by
\L\'od\'z University Grant No.\ 623} \\
Department of Theoretical Physics II \\
University of \L \'od\'z \\
Pomorska 149/153, 90-236 \L \'od\'z, Poland}
\date{}
\maketitle
\setcounter{section}{0}
\setcounter{page}{1}
\begin{abstract}
All deformations of two dimensional centrally extended
Galilei group are classified. The corresponding quantum Lie
algebras are found.
\end{abstract}

\newpage
\section{Introduction}
In this paper we continue the study of the problem of deformations
of spacetime symmetry groups. In recent paper \cite{1} we
classified all nonequivalent Lie-Poisson structures for centrally
extended two-dimensional Galilei group. In the present paper
we quantize these structures by replacing the Poisson brackets
by commutators (this method has been used in Ref.\cite{2} to obtain the
deformations of two-dimensional Galilei group and in Ref.\cite{3}
to construct quantum group dual to the deformed $e_q(2)$).
As a result we
obtain Hopf algebras and find by duality relations the corresponding quantum
Lie algebras.

\section{Deformations of two-dimensional centrally extended Galilei group }
The classical centrally extended Galilei group is defined as a set
of elements
\be
g = (m, \tau , v, a)
\label{r1}
\ee
where $\tau $ is time translation, $a$ and $v$ are space translation
and Galilean boost, respectively, subject to the following multiplication law:
\be
{g}^\prime
{g} =
(m^\prime + m -{1\over 2} v^{\prime 2} \tau -a v^\prime ,
\tau ^\prime  + \tau , v^\prime  + v, a^\prime  + a + \tau v^\prime )
\label{r2}
\ee
The multiplication rule defines coproducts of $m, v, \tau $
and $a$:
\be
\ba{l}
\triangle (m) = m\otimes I + I\otimes m -
{1\over 2}{v^2}\otimes \tau  - v\otimes a \\
\triangle (\tau ) = \tau \otimes I + I\otimes \tau \\
\triangle (a) = a \otimes I + I\otimes a + v \otimes \tau \\
\triangle (v) = v \otimes I + I\otimes v\\
\label{r3}
\ea
\ee
The antipode and counit can be also read off from Eq.(\ref{r2})
\be
\ba{l}
S(m) = -m + {1\over 2} v^{2} \tau -a v \nn\\
S(\tau ) = - \tau \nn\\
S(a) = - a + v \tau \\
S(v) = -v,\nn
\ea
\label{r4}
\ee
\be
\ba{l}
\varepsilon (m) = 0,\;
\varepsilon (\tau ) = 0,\;
\varepsilon (a) = 0,\;
\varepsilon (v) = 0
\label{r5}
\ea
\ee
Galilei group is a real Lie group. The reality condition can be
described by introducing the following $^*$-structure:
\be
m^* = m,\;
\tau ^* = \tau,\;
a^* = a,\;
v^* =v
\label{r6}
\ee
The starting point to obtain two-dimensional quantum Galilei group are the
Lie-Poisson structures on it. In Ref.\cite{1} we have found all relevant
Poisson structures.

In order to obtain the corresponding quantum group we make a replacement
\be
\{ , \} = {1\over {i\kappa }}[ , ]
\label{r7}
\ee
Where $\kappa $ is an arbitrary parameter.
The quantization procedure applied to the Lie-Poisson structures
classified in Ref.\cite{1} yield the noncommutative structures
listed in Table 1.
{\small
\begin{center}
Table 1\\[2,5mm]
{\scriptsize
\begin{tabular}{|c|c|c|c|c|c|c|}
\hline
 &  [ $v$, $a$]& [ $v$, $m$]& [$\tau$, $a$]&
[ $\tau$, $m$]&[ $a$, $m$] &\\ \hline \hline
1&  & &$i{1\over \chi }v$&$-{{1\over 2}i{1\over \chi }{v^2}}$&$
-{{1\over 6}i{1\over \chi }{v^3}}$&${1\over \chi} = {{{\tau _0}^2}\over
\kappa }$ \\  \hline
2&  &$i{1\over \chi }v$& & &$i{1\over \chi }a$&${1\over \chi } ={
{{{{v_0}^2}{\tau _0}} \over \kappa }} $ \\ \hline
3&  & & &$i{1\over \chi }v$&${1\over 2}i{1\over \chi }{v^2}$&${1\over \chi} =
{{{{{v_0}}{{\tau _0}^2}} \over \kappa }} $  \\
\hline
4&  & & & &$-i{1\over \chi }\tau$ &${1\over \chi } ={
{{{{v_0}^3}{\tau _0}} \over \kappa }} $ \\ \hline
5&  & & & &$i{1\over \chi }v$ &${1\over \chi } = {
{{{{v_0}^2}{\tau _0}^2} \over \kappa }} $\\ \hline
6&  &$i{1\over {\chi _1}}v$& & &$i{1\over {\chi _1}}a + i{1\over {\chi _2}}v$
&${1\over \chi_1 } ={
{{{{v_0}^2}{\tau _0}} \over \kappa }}, 
{1\over \chi_2 } = {
{{{{v_0}^2}{\tau _0}^2} \over \kappa }} $ \\ \hline
7 & &$i{1\over {\chi _1}}v$&$i{1\over {\chi _2}}v$&$-{{1\over 2}
i{1\over {\chi _2}}
{v^2}}$&$i{{1\over {\chi _1}}a}-{{1\over 6}i{1\over {\chi _2}}{v^3}}$
&${1\over \chi_1 } ={
{{{{v_0}^2}{\tau _0}} \over \kappa }}, 
{1\over \chi_2 } = {
{{\tau _0}^2} \over \kappa }$  \\ \hline
8&  & & &$i{1\over {\chi _1}}v$&${1\over 2}i{1\over {\chi _1}}{v^2}-
i{1\over {\chi _2}}\tau$  &${1\over \chi_1 } ={
{{{{v_0}}{\tau _0}^2} \over \kappa }}, 
{1\over \chi_2 } ={
{{{{v_0}^3}{\tau _0}} \over \kappa }} $   \\ \hline
9& $i{1\over {\chi _1}}v$&$-{{1\over 2}i{1\over {\chi _1}}{v^2}}$& &
$i{1\over {\chi _2}}v$&$i{1\over {\chi _3}}v-
i{1\over {\chi _1}}m+
{1\over 2}i{1\over {\chi _2}}{v^2}$
&${1\over \chi_1 } ={
{{{{v_0}}{\tau _0}} \over \kappa }}, {1\over \chi_2 } ={
{{{{v_0}}{\tau _0}^2} \over \kappa }}, 
{1\over \chi_3 } ={
{{\varepsilon {{v_0}^2}{\tau _0}^2} \over \kappa }} $   \\ \hline
10&  & &$i{1\over {\chi _1}}v$&$-{{1\over 2}i{1\over {\chi _1}}
{v^2}}$&$-i{{1\over {\chi _2}}\tau}-{1\over 6}i{1\over {\chi _1}}{v^3}$
&${1\over \chi_1 } = {
{{\tau _0}^2} \over \kappa }, {1\over \chi_2 } ={
{{{{v_0}^3}{\tau _0}} \over \kappa }} $  \\
\hline
11&  & & &$-i{1\over {\chi _1}}\tau$&$-i{{1\over {\chi _2}}\tau}
-i{{1\over {\chi _1}}a}$   &${1\over \chi_1} ={
{{{{v_0}^2}{\tau _0}} \over \kappa }}, {1\over \chi_2 } ={
{{{{v_0}^3}{\tau _0}} \over \kappa }} $   \\ \hline
12&  &$i{1\over {\chi _1}}v$& &$-i{1\over {\chi _2}}\tau$
&$i({1\over {\chi _1}} - {1\over {\chi_2}})a$ &${1\over \chi_1 } ={
{{{{v_0}^2}{\tau _0}} \over \kappa }}, {1\over \chi_2 } ={
{{\varepsilon {{v_0}^2}{\tau _0}} \over \kappa }} $   \\ \hline
13&  &$i{1\over {\chi _1}}v$& &$i{{1\over {\chi _1}}\tau}+
i{{1\over {\chi _2}}v}$
&$2i{1\over {\chi _1}}a+{1\over 2}i{1\over {\chi _2}}{v^2}$
&${1\over \chi_1 } ={
{{{{v_0}^2}{\tau _0}} \over \kappa }}, {1\over \chi_2 } ={
{{{{v_0}}{\tau _0}^2} \over \kappa }} $  \\ \hline
14&  &$-i{1\over {\chi _1}}\tau$& &$-i{1\over {\chi _2}}\tau$
&$-i{1\over {\chi _2}}a-{1\over 2}i{1\over {\chi _1}}{{\tau}^2}+
i{1\over {\chi _3}}v$
&${1\over \chi_1 } ={
{{{{v_0}^3}} \over \kappa }}, {1\over \chi_3 } = {
{{{{v_0}^2}{\tau _0}^2} \over \kappa }}, 
{1\over \chi_2 } ={
{{\varepsilon {{v_0}^2}{\tau _0}} \over \kappa }} $  \\ \hline
15&  &$-i{1\over {\chi _1}}\tau$& &$-i{1\over {\chi _2}}\tau+
i{1\over {\chi _3}}v$
&$-i{1\over {\chi _2}}a-{1\over 2}i{1\over {\chi _1}}{{\tau}^2}+
{1\over 2}
i{1\over {\chi _3}}{{v}^2}$
&${1\over \chi_1 } ={
{{{{v_0}^3}} \over \kappa }}, {1\over \chi_3 } = {
{{{{v_0}}{\tau _0}^2} \over \kappa }}, 
{1\over \chi_2 } ={
{{\varepsilon {{v_0}^2}{\tau _0}} \over \kappa }} $  \\ \hline
16&  & &$-i{1\over {\chi }}\tau$&$i{1\over {\chi }}a$&$
-i{1\over {\chi }}m$&$
{1\over \chi } = {{{v_0}{\tau _0}}\over \kappa }$ \\ \hline
\end{tabular}
}
\end{center}
}
and, for all cases:
\be
[v, \tau ] = 0
\label{r8}
\ee
We have checked that the comutators listed
above satisfy the Jacobi identities.
We have also verified that the equations
(\ref{r3}), (\ref{r4}), (\ref{r5}), (\ref{r6}) together with Table 1 and
(\ref{r8}) define a $ ^*$-Hopf algebras which provide the deformations
of two-dimensional centrally extended Galilei group.
\section{The quantum Lie algebra and coalgebra}
In this section we find all quantum Lie algebras corresponding to
the quantum Galilei groups defined in Sec.2.
To this end we use the Hopf algebra duality rules.
On the classsical level the generators of Lie algebra of Galilei group
can be defined by the following global parametrization of group
element
\be
g = e^{imM} e^{-i\tau H} e^{iaP} e^{ivK}
\label{r9}
\ee
We adopt here this definition as well as classical duality relation
\be
<\varphi ,X> = -i{d\over dt} \varphi (e^{itX}) {\mid_{t=o}}
\label{r10}
\ee
The group algebra is generated by the set of elements of the form
\be
\varphi  = {m^\alpha }{\tau ^\beta }{a^\gamma }{v^\sigma}
\label{r11}
\ee
where ${\alpha ,\beta ,\gamma ,\sigma} \geq 1 $.
By applying the classical duality rules we obtain
\be
\ba{l}
<M, {m^\alpha }{\tau ^\beta }{a^\gamma }{v^\sigma}> = - i \delta _{1\alpha }
\delta _{0\beta } \delta _{0\gamma }\delta _{0\sigma }\\
<H, {m^\alpha }{\tau ^\beta }{a^\gamma }{v^\sigma}> = i \delta _{0\alpha }
\delta _{1\beta } \delta _{0\gamma }\delta _{0\sigma } \\
<P, {m^\alpha }{\tau ^\beta }{a^\gamma }{v^\sigma}> = - i \delta _{0\alpha }
\delta _{0\beta } \delta _{1\gamma }\delta _{0\sigma }\\
<K, {m^\alpha }{\tau ^\beta }{a^\gamma }{v^\sigma}> = - i \delta _{0\alpha }
\delta _{0\beta } \delta _{0\gamma }\delta _{1\sigma }.\\
\label{r12}
\ea
\ee
In order to define dual algebra structures we use
the duality relation:
\be
\ba{l}
<\varphi , XY> = <{\triangle \varphi}, X \otimes Y>\\
<\varphi  \psi, X> = <\varphi  \otimes \psi ,\triangle (X)>
\label{r13}
\ea
\ee

The ${^*}$-structure on dual Hopf algebra can be defined by the formula
\be
<X^*, \varphi > = <X, {S^{-1}}(\varphi ^*)>
\label{r14}
\ee
provided the following relation holds
\be
S^{-1}(X) = [S(X^*)]^*
\label{r15}
\ee
In order to take care about all elements of the form (\ref{r11})
let us introduce the generating function depending on
four real parameters $\mu ,\nu , \varrho, \kappa $
\be
\varphi (\mu, \nu, \rho, \kappa ) = <\varphi , e^{\mu m} e^{\nu \tau }
e^{\rho a} e^{\kappa  v}>
\label{r16}
\ee
Then we have
\be
\ba{l}
<M,e^{\mu m} e^{\nu \tau }
e^{\rho a} e^{\kappa  v} > =  - i\mu\\
<H,e^{\mu m} e^{\nu \tau }
e^{\rho a} e^{\kappa  v} > = i\nu \\
<P,e^{\mu m} e^{\nu \tau }
e^{\rho a} e^{\kappa  v} > = -i\varrho \\
<K,e^{\mu m} e^{\nu \tau }
e^{\rho a} e^{\kappa  v} > = - i\kappa \\
\label{r17}
\ea
\ee

The duality rules (\ref{r13}) and long and tedious
calculations lead us to the following structure of quantum
Lie algebras
\be
\ba{l}
1.\\
\;\;{[K, P] = iM} \\
\;\;{[K, H] = iP} \\
\\
\;\;\Delta M = I \otimes M + M \otimes I\\
\;\;\Delta H = I \otimes H + H \otimes I\\
\;\;\Delta P = I \otimes P + P \otimes I\\
\;\;\Delta K = I \otimes K + K \otimes I - {1\over {\chi} }P \otimes H\\
\\
\;\;S(M) =  - M\\
\;\;S(H) = - H\\
\;\;S(P) =  -P\\
\;\;S(K) = - K - {{1\over \chi}PH}
\label{r18}
\ea
\ee
\be
\ba{l}
2.\\
\;\;{[K, H] = iP}\\
\;\;{[K, P] = {1\over 2}i\chi {(1-
{e^{-2{1\over \chi }M}})}} \\
\\
\;\;\Delta M = I \otimes M + M \otimes I\\
\;\;\Delta H = I \otimes H + H \otimes I\\
\;\;\Delta P = I \otimes P + P \otimes e^{-{1\over {\chi} }M}\\
\;\;\Delta K = I \otimes K + K \otimes e^{-{1\over {\chi} }M}\\
\\
\;\;S(M) =  - M\\
\;\;S(H) = - H\\
\;\;S(P) =  -Pe^{{1\over \chi}M}\\
\;\;S(K) = - Ke^{{1\over \chi}M}
\label{r19}
\ea
\ee
\be
\ba{l}
3.\\
\;\;{[K, P] = iM} \\
\;\;{[K, H] = iP} \\
\\
\;\;\Delta M = I \otimes M + M \otimes I\\
\;\;\Delta H = I \otimes H + H \otimes I\\
\;\;\Delta P = I \otimes P + P \otimes I\\
\;\;\Delta K = I \otimes K + K \otimes I + {1\over {\chi} }H \otimes M\\
\\
\;\;S(M) =  - M\\
\;\;S(H) = - H\\
\;\;S(P) =  -P\\
\;\;S(K) = - K - {1\over \chi}HM 
\label{r20}
\ea
\ee
\be
\ba{l}
4.\\
\;\;{[K, P] = iM}\\
\;\;{[K, H] = iP - {{1\over 2}{i{1\over \chi }{M^2}}} } \\
\\
\;\;\Delta M = I\otimes M + M \otimes I\\
\;\;\Delta H = I \otimes H + H \otimes I - {1\over {\chi} }P \otimes M\\
\;\;\Delta P = I \otimes P + P \otimes I\\
\;\;\Delta K = I \otimes K + K \otimes I\\
\\
\;\;S(M) =  - M\\
\;\;S(H) = - H - {1\over \chi}PM\\
\;\;S(P) =  -P\\
\;\;S(K) = - K
\label{r21}
\ea
\ee
\be
\ba{l}
5.\\
\;\;{[K, P] = iM}\\
\;\;{[K, H] = iP}\\
\\
\;\;\Delta M = I \otimes M + M \otimes I\\
\;\;\Delta H = I \otimes H + H \otimes I\\
\;\;\Delta P = I \otimes P + P \otimes I\\
\;\;\Delta K = I \otimes K + K \otimes I -
{1\over {\chi} }P \otimes M\\
\\
\;\;S(M) =  - M\\
\;\;S(H) = - H\\
\;\;S(P) =  -P\\
\;\;S(K) = - K - {1\over \chi}PM
\label{r22}
\ea
\ee
\be
\ba{l}
6.\\
\;\;{[K, H] = iP}\\
\;\;{[K, P] = {1\over 2}i {\chi_1} {(1-
{e^{-2{1\over {\chi_1} }M}})}}  \\
\\
\;\;\;\Delta M = I \otimes M + M \otimes I\\
\;\;\;\Delta H = I \otimes H + H \otimes I\\
\;\;\;\Delta P = I \otimes P + P \otimes e^{-{1\over {\chi_1} }M}\\
\;\;\;\Delta K = I \otimes K + K \otimes e^{-{1\over {\chi_1} }M}
- {1\over {\chi_2} }P \otimes Me^{-{1\over {\chi_1} }M}\\
\\
\;\;S(M) =  - M\\
\;\;S(H) = - H\\
\;\;S(P) =  -Pe^{{1\over \chi_1}M}\\
\;\;S(K) = - Ke^{{1\over \chi_1}M} - {1\over \chi_2}PM
e^{{1\over \chi_1}M}
\label{r23}
\ea
\ee
\be
\ba{l}
7.\\
\;\;{[K, H] = iP}\\
\;\;{[K, P] = {1\over 2} i {\chi_1} {(1-
{e^{-2{1\over {\chi_1} }M}})} } \\
\\
\;\;\;\Delta M = I \otimes M + M \otimes I\\
\;\;\;\Delta H = I \otimes H + H \otimes I\\
\;\;\;\Delta P = I \otimes P + P \otimes e^{-{1\over {\chi_1} }M}\\
\;\;\;\Delta K = I \otimes K + K \otimes e^{-{1\over {\chi_1} }M}
- {1\over {\chi_2} }P \otimes He^{-{1\over {\chi_1} }M}\\
\\
\;\;S(M) =  - M\\
\;\;S(H) = - H\\
\;\;S(P) =  -Pe^{{1\over \chi_1}M}\\
\;\;S(K) = - Ke^{{1\over \chi_1}M} - {1\over \chi_2}PH
e^{{1\over \chi_1}M}
\label{r24}
\ea
\ee
\be
\ba{l}
8.\\
\;\;{[K, P] = iM}\\
\;\;{[K, H] = iP - {{1\over 2}{i{1\over {\chi_2} }{M^2}}} } \\
\\
\;\;\Delta M = I \otimes M + M \otimes I\\
\;\;\Delta H = I \otimes H + H \otimes I - {1\over {\chi_2} }P \otimes M\\
\;\;\Delta P = I \otimes P + P \otimes I\\
\;\;\Delta K = I \otimes K + K \otimes I + {1\over {\chi_1} }H \otimes M
- {1\over 2}{1\over {\chi_1} }{1\over {\chi_2} }P \otimes {M^2}\\
\\
\;\;S(M) =  - M\\
\;\;S(H) = - H - {1\over \chi_2}PM\\
\;\;S(P) =  -P\\
\;\;S(K) = - K + {1\over \chi_1}HM + {1\over 2}{1\over \chi_1}
{1\over \chi_2}P{M^2}
\label{r25}
\ea
\ee
\be
\ba{l}
9.\\
\;\;{[K, M]} = {{1\over 2}{i{1\over {\chi_1} }{M^2}} } \\
\;\;{[K, P] = iM}\\
\;\;{[K, H] = - i{\chi_1}(1 - e^{{1\over {\chi_1} }P})}\\
\\
\;\;\Delta M = {e^{{1\over {\chi_1} }P}} \otimes M + M \otimes I\\
\;\;\Delta H = I \otimes H + H \otimes I\\
\;\;\Delta P = I \otimes P + P \otimes I\\
\;\;\Delta K = {e^{{1\over {\chi_1} }P}} \otimes K + K \otimes I
+ {1\over {\chi_2} }H{e^{{1\over {\chi_1} }P}} \otimes M
- {1\over {\chi_3} }P{e^{{1\over {\chi_1} }P}} \otimes M\\
\\
\;\;S(M) = - Me^{- {1\over \chi_1}P}\\
\;\;S(H) = - H\\
\;\;S(P) = - P\\
\;\;S(K) = - Ke^{- {1\over \chi_1}P} + {1\over \chi_2}HM
e^{-{1\over \chi_1}P} - {1\over \chi_3}PM
e^{-{1\over \chi_1}P}
\label{r26}
\ea
\ee
\be
\ba{l}
10.\\
\;\;{[K, P] = iM}\\
\;\;{[K, H] = iP - {{1\over 2}{i{{1\over {\chi _2}} }{M^2}} }} \\
\\
\;\;\Delta M = I \otimes M + M \otimes I\\
\;\;\Delta H = I \otimes H + H \otimes I - {1\over {\chi_2}}P \otimes M\\
\;\;\Delta P = I \otimes P + P \otimes I\\
\;\;\Delta K = I \otimes K + K \otimes I
+ {{1\over 2}{1\over {\chi_1} }{1\over {\chi_2} }{P^2}} \otimes M
- {1\over {\chi_1} }P \otimes H\\
\\
\;\;S(M) =  - M\\
\;\;S(H) = - H - {1\over \chi_2}PM\\
\;\;S(P) =  -P\\
\;\;S(K) = - K - {1\over \chi_1}PM - {1\over 2}{1\over \chi_1}
{1\over \chi_2}M{P^2}
\label{r27}
\ea
\ee
\be
\ba{l}
11.\\
\;\;{[K, H] = iP - i{{\chi _1}\over {\chi _2}}M
{e^{{1\over {\chi _1}}M}} + i{{\chi _1}^2 \over {\chi _2}}
(e^{{1\over {\chi _1}}M} - 1)}\\
\;\;{[K, P] = -i{\chi_1}
(e^{{1\over {\chi _1}}M} - 1)}\\
\\
\;\;\;\Delta M = I \otimes M + M \otimes I\\
\;\;\;\Delta H = I \otimes H + H \otimes e^{{1\over {\chi_1} }M} -
{1\over {\chi_2} }P \otimes Me^{{1\over {\chi_1} }M}\\
\;\;\;\Delta P = I \otimes P + P \otimes e^{{1\over {\chi_1} }M}\\
\;\;\;\Delta K = I \otimes K + K \otimes I\\
\\
\;\;S(M) =  - M\\
\;\;S(H) = - He^{-{1\over \chi_1}M} - {1\over \chi_2}PM
e^{-{{1\over \chi_1}M}}\\
\;\;S(P) =  -Pe^{-{1\over \chi_1}M}\\
\;\;S(K) = - K
\label{r28}
\ea
\ee
\be
\ba{l}
12.\\
\;\;{[K, H] = iP}\\
\;\;{[K, P] = {-{i\over {({2{1\over \chi_1} - {1\over \chi_2}})}}
{(1 - {e^{-({2{1\over \chi_1} - {1\over \chi_2}})M}}})} }\\
\\
\;\;\Delta M = I \otimes M + M \otimes I\\
\;\;\Delta H = I \otimes H + H \otimes {e^{{1\over {\chi_2} }M}}\\
\;\;\Delta P = I \otimes P + P \otimes {e^{({1\over {\chi_2} } -
{1\over {\chi_1} })M}}\\
\;\;\Delta K = I \otimes K + K \otimes {e^{-{1\over {\chi_1} }M}}\\
\\
\;\;S(M) =  - M\\
\;\;S(H) = - He^{-{1\over \chi_2}M}\\
\;\;S(P) =  -Pe^{({1\over \chi_1} - {1\over \chi_2})M}\\
\;\;S(K) = - K e^{{1\over \chi_1}M}
\label{29}
\ea
\ee
\be
\ba{l}
13.\\
\;\;{[K, H] = iP}\\
\;\;[K, P] = {1\over 3}i{\chi _1}{(1-
{e^{-3{1\over {\chi _1}}M}})}  \\
\\
\;\;\Delta M = I \otimes M + M \otimes I\\
\;\;\Delta H = I \otimes H + H \otimes {e^{-{1\over {\chi_1} }M}}\\
\;\;\Delta P = I \otimes P + P \otimes {e^{-2{1\over {\chi_1} }M}}\\
\;\;\Delta K = I \otimes K + K \otimes {e^{-{{1\over {\chi_1} }M}}}
+ {1\over {\chi_2} }H \otimes M{e^{-{1\over {\chi_1} }M}}\\
\\
\;\;S(M) =  - M\\
\;\;S(H) = - He^{{1\over \chi_1}M}\\
\;\;S(P) =  -Pe^{2{1\over \chi_1}M}\\
\;\;S(K) = - K e^{{1\over \chi_1}M} + {1\over \chi_2}HM
e^{{1\over \chi_1}M}
\label{r30}
\ea
\ee
\be
\ba{l}
14.\\
\;\;{[K, H] = iP + i{{\chi_2}^2 \over {{\chi_1}{\chi_3}}}M
{e^{{1\over {\chi _2}}M}} - {1\over 2}i{{{\chi_2}^3 \over {{\chi_1}{\chi_3}}}}
(e^{2{1\over {\chi _2}}M}} - 1)\\
\;\;{[K, P] = i{\chi_2}
(e^{{1\over {\chi _2}}M}} - 1)\\
\;\;{[H, P] = - {1\over 2}i{{{\chi_2}^2 \over {{\chi_1}}}}
(e^{2{1\over {\chi _2}}M}} - 1) +
i{{{\chi_2}^2 \over {{\chi_1}}}}
(e^{{1\over {\chi _2}}M} - 1) \\
\\
\;\;\;\Delta M = I \otimes M + M \otimes I\\
\;\;\;\Delta H = I \otimes H + H \otimes e^{{1\over {\chi_2} }M} -
{{\chi_2}\over {\chi_1} }K \otimes
e^{{1\over {\chi_2} }M} +
{{\chi_2}\over {\chi_1} }K \otimes I \\
\;\;+ {{{\chi_2}\over {{\chi_1}{\chi_3}} }P} \otimes
{Me^{{1\over {\chi_2} }M}} + {{{\chi_2}^2}
\over {{\chi_1}{\chi_3}} }P
\otimes I - {{{\chi_2}^2}\over {{\chi_1}{\chi_3}} }P \otimes
e^{{1\over {\chi_2} }M}\\
\;\;\;\Delta P = I \otimes P + P \otimes e^{{1\over {\chi_2} }M}\\
\;\;\;\Delta K = I \otimes K + K \otimes I + {{{\chi_2}\over {\chi_3} }P
\otimes I} - {{{\chi_2}\over {\chi_3} }P \otimes
e^{{1\over {\chi_2} }M}}\\
\\
\;\;S(M) =  - M\\
\;\;S(H) = - He^{-{1\over \chi_2}M} - {\chi_2 \over \chi_1}K
+ {\chi_2 \over \chi_1}
Ke^{-{1\over \chi_2}M} + {{\chi_2}^2 \over {\chi_1 \chi_3}}
Pe^{-{1\over \chi_2}M} \\
\;\;-
{{\chi_2}^2 \over {\chi_1 \chi_3}}P + {{\chi_2} \over {\chi_1 \chi_3}}
PMe^{-{1\over \chi_2}M}\\
\;\;S(P) =  -Pe^{-{1\over \chi_2}M}\\
\;\;S(K) = - K + {{\chi_2} \over { \chi_3}}
Pe^{-{1\over \chi_2}M} - {{\chi_2} \over { \chi_3}}P 
\label{31}
\ea
\ee
\be
\ba{l}
15.\\
\;\;{[K, H] = iP}\\
\;\;{[K, P] =  {
i
\left(
- {1\over 2}{1\over \chi_2} + {1\over 2}
\sqrt{ {1\over {\chi_2}^2} - 4{1\over {\chi_1}}
{1\over {\chi_3}}}
\right)
\over {\sqrt{ {1\over {\chi_2}^2} - 4{1\over {\chi_1}}
{1\over {\chi_3}}}
\left(
{3\over 2}{1\over \chi_2} + {1\over 2}
\sqrt{ {1\over {\chi_2}^2} - 4{1\over {\chi_1}}
{1\over {\chi_3}}}
\right)  }
}
\left(
e^{({3\over 2}{1\over \chi_2} + {1\over 2}
\sqrt{ {1\over {\chi_2}^2} - 4{1\over {\chi_1}}
{1\over {\chi_3}}})M} - 1
\right)}\\
\;\;+ {
i
\left(
 {1\over 2}{1\over \chi_2} + {1\over 2}
\sqrt{ {1\over {\chi_2}^2} - 4{1\over {\chi_1}}
{1\over {\chi_3}}}
\right)
\over {\sqrt{ {1\over {\chi_2}^2} - 4{1\over {\chi_1}}
{1\over {\chi_3}}}
\left(
{3\over 2}{1\over \chi_2} - {1\over 2}
\sqrt{ {1\over {\chi_2}^2} - 4{1\over {\chi_1}}
{1\over {\chi_3}}}
\right)  }
}
\left(
e^{({3\over 2}{1\over \chi_2} - {1\over 2}
\sqrt{ {1\over {\chi_2}^2} - 4{1\over {\chi_1}}
{1\over {\chi_3}}})M} - 1
\right)\\
\;\;{[H, P] = {
i\over {\chi_1 \sqrt{ {1\over {\chi_2}^2} - 4{1\over {\chi_1}}
{1\over {\chi_3}}}
\left(
{3\over 2}{1\over \chi_2} - {1\over 2}
\sqrt{ {1\over {\chi_2}^2} - 4{1\over {\chi_1}}
{1\over {\chi_3}}}
\right)  }
}
\left(
e^{({3\over 2}{1\over \chi_2} - {1\over 2}
\sqrt{ {1\over {\chi_2}^2} - 4{1\over {\chi_1}}
{1\over {\chi_3}}})M} - 1
\right)}\\
\;\;- {
i\over {\chi_1 \sqrt{ {1\over {\chi_2}^2} - 4{1\over {\chi_1}}
{1\over {\chi_3}}}
\left(
{3\over 2}{1\over \chi_2} + {1\over 2}
\sqrt{ {1\over {\chi_2}^2} - 4{1\over {\chi_1}}
{1\over {\chi_3}}}
\right)  }
}
\left(
e^{({3\over 2}{1\over \chi_2} + {1\over 2}
\sqrt{ {1\over {\chi_2}^2} - 4{1\over {\chi_1}}
{1\over {\chi_3}}})M} - 1
\right)\\
\\
\;\;\;\Delta M = I \otimes M + M \otimes I\\
\;\;\;\Delta H = I \otimes H + H \otimes e^{{1\over 2}{1\over {\chi_2}}M}
[\cosh({-{1\over 2}\sqrt{ {1\over {\chi_2}^2} - 4{1\over {\chi_1}}
{1\over {\chi_3}}}M})
\\
\;\;-  {1\over {{\chi_2}\sqrt{ {1\over {\chi_2}^2} - 4{1\over {\chi_1}}
{1\over {\chi_3}}}}}
\sinh ({-{1\over 2}\sqrt{ {1\over {\chi_2}^2} - 4{1\over {\chi_1}}
{1\over {\chi_3}}}M})]\\
\;\; +
2 {1\over {{\chi_1}\sqrt{ {1\over {\chi_2}^2} - 4{1\over {\chi_1}}
{1\over {\chi_3}}}}} K \otimes
e^{{1\over 2}{1\over {\chi_2}}M}
\sinh({-{1\over 2}\sqrt{ {1\over {\chi_2}^2} - 4{1\over {\chi_1}}
{1\over {\chi_3}}}M})
\label{r32}
\ea
\ee
\be
\ba{l}
\;\;\;\Delta P = I \otimes P + P \otimes e^{{1\over {\chi_2}}M}\\
\;\;\;\Delta K = I \otimes K + K \otimes {e^{{{1\over 2}{1\over {\chi_2}}M}}
{[\cosh({-{1\over 2}\sqrt{ {1\over {\chi_2}^2} - 4{1\over {\chi_1}}
{1\over {\chi_3}}}M})}}
\\
\;\;+  {1\over {{\chi_2}\sqrt{ {1\over {\chi_2}^2} - 4{1\over {\chi_1}}
{1\over {\chi_3}}}}}
\sinh ({-{1\over 2}\sqrt{ {1\over {\chi_2}^2} - 4{1\over {\chi_1}}
{1\over {\chi_3}}}M})] \\
\;\; -
2 {1\over {{\chi_3}\sqrt{ {1\over {\chi_2}^2} - 4{1\over {\chi_1}}
{1\over {\chi_3}}}}} H \otimes
e^{{1\over 2}{1\over {\chi_2}}M}
\sinh({-{1\over 2}\sqrt{ {1\over {\chi_2}^2} - 4{1\over {\chi_1}}
{1\over {\chi_3}}}M})\\
\\
\;\;S(M) =  - M\\
\;\;S(H) = - He^{-{1\over 2}{1\over \chi_2}M}[\cosh
({{1\over 2}\sqrt{ {1\over {\chi_2}^2} - 4{1\over {\chi_1}}
{1\over {\chi_3}}}M})\\
\;\;
-{1\over {{\chi_2}\sqrt{ {1\over {\chi_2}^2} - 4{1\over {\chi_1}}
{1\over {\chi_3}}}}}\sinh
({{1\over 2}\sqrt{ {1\over {\chi_2}^2} - 4{1\over {\chi_1}}
{1\over {\chi_3}}}M})] \\
\;\;-
{2{1\over {{\chi_1}\sqrt{ {1\over {\chi_2}^2} - 4{1\over {\chi_1}}
{1\over {\chi_3}}}}}}K
e^{-{1\over 2}{1\over \chi_2}M}
\sinh ({{1\over 2}\sqrt{ {1\over {\chi_2}^2} - 4{1\over {\chi_1}}
{1\over {\chi_3}}}M})\\
\;\;S(P) =  -Pe^{-{1\over \chi_2}M}\\
\;\;S(K) = - Ke^{-{1\over 2}{1\over \chi_2}M}[\cosh
({{1\over 2}\sqrt{ {1\over {\chi_2}^2} - 4{1\over {\chi_1}}
{1\over {\chi_3}}}M})\\
\;\;
+ {1\over {{\chi_2}\sqrt{ {1\over {\chi_2}^2} - 4{1\over {\chi_1}}
{1\over {\chi_3}}}}}\sinh
({{1\over 2}\sqrt{ {1\over {\chi_2}^2} - 4{1\over {\chi_1}}
{1\over {\chi_3}}}M})] \\
\;\;+
{2{1\over {{\chi_3}\sqrt{ {1\over {\chi_2}^2} - 4{1\over {\chi_1}}
{1\over {\chi_3}}}}}}H
e^{-{1\over 2}{1\over \chi_2}M}
\sinh ({{1\over 2}\sqrt{ {1\over {\chi_2}^2} - 4{1\over {\chi_1}}
{1\over {\chi_3}}}M})
\label{r33}
\ea
\ee
\be
\ba{l}
16.\\
\;\;{[K, M] = {1\over 2}i{1\over \chi}{M^2}}\\
\;\;{[K, P] = iM }\\
\;\;{[K, H] = i{\chi}(1 - e^{-{1\over \chi}P}) - i
{1\over \chi}MH}\\
\\
\;\;\Delta M = M \otimes I + {{e^{{1\over {\chi}} P} \otimes M}\over {1-
{1\over 2}{1\over {\chi^2}}He^{{1\over {\chi}} P} \otimes M}}\\
\;\;\Delta H = e^{-{1\over {\chi}} P} \otimes H + H \otimes I
- {{1\over {\chi}^2} } H \otimes MH + {1\over 4}{{1\over {\chi}^4} }
{H^2}e^{{1\over {\chi}} P} \otimes {M^2}H \\
\;\;- {1\over 2}{{1\over {\chi}^2} }
{H^2}e^{{1\over {\chi}} P} \otimes M\\
\;\;\Delta P = I \otimes P + P \otimes I - {2\chi }\ln
{(1 - {1\over 2}{1\over {\chi}^2} {H}e^{{1\over {\chi}} P} \otimes M)}\\
\;\;\Delta K = I \otimes K + K \otimes I\\
\\
\;\;S(M) =  {{- Me^{-{1\over \chi}P}}\over {1 - {1\over 2}
{1\over \chi^2}MH}}\\
\;\;S(H) = - He^{{1\over \chi}P} + {1\over 2}{1\over \chi^2}M{H^2}
e^{{1\over \chi}P}\\
\;\;S(P) =  -P + {2\chi }\ln \left(
{1\over {1 - {1\over 2}{{1\over \chi^2}}HM}}
\right)\\
\;\;S(K) = - K
\label{r34}
\ea
\ee
The counits for all cases take the form
\be
\ba{l}
\varepsilon (M) = 0,\;
\varepsilon (H ) = 0,\;
\varepsilon (P) = 0,\;
\varepsilon (K) = 0
\label{r35}
\ea
\ee
\section{Conclusions}
We obtained a number of, in general multiparameter, deformations of
two-dimensional centrally extended Galilei group. They are described
in Table 1 together with Eqs. (\ref{r3}), (\ref{r4}), (\ref{r5}),
(\ref{r6}) and (\ref{r8}). The corresponding quantum Lie  algebras
have been also found. They are listed in Eqs.(\ref{r18}) - (\ref{r35}).
Some of the resulting structures seem to be quite interesting and
deserve more detailed study which will be the subject of subsequent
publications.
\section{Acknowledgment}
The Author acknowledges Prof. S.~Giller and Prof.~P.~Kosi\'nski
for a careful reading of the
manuscript and many helpful suggestion. Special thanks are also due to
Prof. P. Ma\'slanka and Mr.~E.~Kowalczyk for valuable discussion.
\section{Appendix}
In this appendix, in order to illustrate the procedures used,
we obtain Eqs. (\ref{r26})
for the case 9.

In this case the generating function reads (cf. remark after
Eq. (\ref{r16}))
\[
\varphi = e^{\mu m} e^{\nu \tau }
e^{\kappa  v} e^{\rho a}
\nn
\]
In order to determine coalgebra sector we use the duality
relation\\
\[
<\varphi  \psi, X> = <\varphi  \otimes \psi ,\triangle (X)>
\]
To this end, we calculate:
{\large
\[
e^{\mu  ^\p m} e^{\nu ^\p \tau } e^{\kappa  ^\p v} e^{\rho ^\p a}
e^{\mu m} e^{\nu \tau } e^{\kappa  v} e^{\rho a} =
e^{(\mu ^\p + \mu e^{-i{1\over {\chi_1}}\varrho ^\p})m}
e^{(\nu ^\p + \nu)\tau }
\]
\[
e^{
\left(
\kappa  e^{-i{1\over {\chi_1}}\varrho ^\p} + {\kappa ^\p}
{1\over {{1\over 2}i{1\over {\chi_1}}\mu
e^{-i{1\over {\chi_1}}\varrho ^\p}v}
+ 1}
\right) v + (2{{{\chi_1}\over {\chi_2}}}{\nu ^\p} + 2
{\chi_1 \over \chi_3}{\rho ^\p})
\ln{
\left(
{1\over 2}i{1\over {\chi_1}}v\mu e^{-i{1\over {\chi_1}}{\rho ^\p}} + 1
\right) }
}
\]
\[
e^{O({v^2}, ...)}
e^{({\rho ^\p} + \rho )a}
\]     }
where $O({v^2}, ...)$ means function depending on $v^2$
and higher powers of $v$.

Taking the first power of $v$ in exponentials we obtain the following result\\
{\large
\[
\varphi \psi =
e^{(\mu ^\p + \mu e^{-i{1\over {\chi_1}}\varrho ^\p})m}
e^{(\nu ^\p + \nu)\tau }
\]
\[
e^{
\left(
\kappa  e^{-i{1\over {\chi_1}}\varrho ^\p} + {\kappa ^\p}
+  {i{1\over {\chi_2}}\mu {\nu ^\p}e^{-i{1\over {\chi_1}}\varrho ^\p}}
+  i{1\over {\chi_3}}\mu {\rho ^\p}e^{-i{1\over {\chi_1}}\varrho ^\p}
\right) v
}
e^{O({v^2}, ...)}
e^{({\rho ^\p} + \rho )a}
\]
}
Therefore we have
\[
<\Delta M, e^{\mu  ^\p m} e^{\nu ^\p \tau } e^{\kappa  ^\p v} e^{\rho ^\p a}
\otimes e^{\mu m} e^{\nu \tau } e^{\kappa  v} e^{\rho a} > =
<M, e^{\mu  ^\p m} e^{\nu ^\p \tau } e^{\kappa  ^\p v} e^{\rho ^\p a}
e^{\mu m} e^{\nu \tau } e^{\kappa  v} e^{\rho a} > =
\]
$
\;\;\; = - i(\mu ^\p + \mu e^{-i{1\over {\chi_1}}\varrho ^\p})
$\\

Using this and the corresponding formulas for $H, P$ and $K$
we arrive at the coalgebra described by Eq.(\ref{r26})\\

Now using the duality relation $<\varphi , XY> =
<{\triangle \varphi}, X \otimes Y>$ we can determine the algebra sector.
We have:\\
{\large
\[
\Delta{\varphi} = e^{\mu (m^\prime + m -{1\over 2} v^{\prime 2} \tau
-a v^\prime )}
e^{\nu (\tau ^\prime  + \tau ) }
e^{\kappa  (v^\prime  + v)}
e^{\rho (a^\prime  + a + \tau v^\prime)}
\nn
\]
}
In order to deal with the above expression we first decompose
{\large
\[
e^{\mu (m^\prime + m -{1\over 2} v^{\prime 2} \tau -a v^\prime )}
\]
}
To this end let
{\large
\[
f(\mu ) =  e^{-\mu {m^\p}}e^{\mu (m^\prime + m -{1\over 2} v^{\prime 2}
\tau -a v^\prime )}.
\nn
\]
}
Then
{\large
\[
\Delta{\varphi} = e^{\mu m^\prime }f
e^{\nu (\tau ^\prime  + \tau ) }
e^{\kappa  (v^\prime  + v)}
e^{\rho (a^\prime  + a + \tau v^\prime)}.
\nn
\]
}
The function $f(\mu )$ obeys $f(0) = 1$, and
\[
\dot{f} = \frac{df}{dx} =
\left(
m - {1\over 2}{{\tau {{v^\p}^2}}\over (1 + {1\over 2}i{1\over {\chi_1}}\mu
{v^\p})^2} - {{a{v^\p}}\over {1 + {1\over 2}i{1\over {\chi_1}}\mu
{v^\p}}}
\right)
f
\]
The above equation cannot be solved in a standard
way due to the fact that the terms appearing on the right
hand side do not commute. Therefore we pass to the "interaction picture"
by letting
{\large
\[
f = e^{mX(v^\p , \mu )}h.
\]
}
Then
{\large
\[
 \dot{f} = m\dot{X}e^{mX}h + e^{mX}\dot{h}
\]
}
and we select $X$ in such a way that the terms containing $m$ cancel. It can
be checked that $X$ should be of the form
\[
X(v^\p , \mu ) = \mu (1 + {1\over 2}i{1\over {\chi_1}}\mu
{v^\p}).
\]
We have now:
{\large
\[
e^{\mu (m^\prime + m -{1\over 2} v^{\prime 2} \tau -a v^\prime )}
= e^{\mu {m^\p}}e^{\mu m(1 + {1\over 2}i{1\over {\chi_1}}\mu
{v^\p})}h
\]
}
Our next aim is  calculate $h$. We do this in the same way by
writing out a relevant differential equation
\[
\dot{h} = -
\left(
{{a{v^\p}}\over {1 + {1\over 2}i{1\over {\chi_1}}\mu
{v^\p}}} + \gamma (\tau , v, v^\p , \mu )
\right)
h
\]
and substituting
{\large
\[
h = e^{\lambda ({v^\p}, \mu )a}
\]
}
in order to deal with noncommutativity of the terms on
right hand side.

Step by step we arrive at following result
{\large
\[
e^{\mu (m^\prime + m -{1\over 2} v^{\prime 2} \tau -a v^\prime )}
=
\]
\[ 
e^{\mu {m^\p}}e^{\mu m(1 + {1\over 2}i{1\over {\chi_1}}\mu
{v^\p})} e^{{{2i}{{\chi_1}}}\ln
(1 + {1\over 2}i{1\over {\chi_1}}\mu
{v^\p})a}
e^{-\int_{0}^{\mu }{\gamma ({v\over {
(1 + {1\over 2}i{1\over {\chi_1}}\mu
{v^\p})^2}}, {v^\p}, \tau , \mu )}d\mu }
\]
}
where funtion $\gamma $ is of third degree in $v, v^\p$.

Finally, $\Delta\varphi $ can be rewritten as
{\large
\[
\Delta\varphi = 
e^{\mu m^\p}e^{\nu \tau ^\p}e^{\mu m}
e^{{1\over 2}i{1\over {\chi_1}}{\mu ^2}{v\p }m}e^{\nu \tau}
e^{\kappa {v^\p}}e^{O({{v^\p}^2}, v, \tau )}
\]
\[
e^{\kappa ve^{i{1\over {\chi_1}}\mu
{v^\p}}}e^{-\mu {v^\p}a}e^{\varrho a}e^{-i\chi
(1 - e^{-i{1\over {\chi_1}}\varrho }){v^\p}\tau }e^{\varrho {a^\p}}
\]
}
Using the above formula, we may calculate
\[
<KM, e^{\mu m} e^{\nu \tau } e^{\kappa  v} e^{\rho a} > =
<K\otimes M, \Delta(e^{\mu m} e^{\nu \tau } e^{\kappa  v} e^{\rho a}) > =
\]
\[
= (-i)(-i)({1\over 2}i{1\over {\chi_1}}{\mu ^2})\\
\]
\[
<MK, e^{\mu m} e^{\nu \tau } e^{\kappa  v} e^{\rho a} > =
<M\otimes K, \Delta(e^{\mu m} e^{\nu \tau } e^{\kappa  v} e^{\rho a}) > = 0
\nn
\]
Therefore:
\[
<[K, M], e^{\mu m} e^{\nu \tau } e^{\kappa  v} e^{\rho a} > =
-{1\over 2}i{1\over {\chi_1}}{\mu ^2}
\]
Other commutators are calculated in a similar way.



\begin{thebibliography}{99}
\newcommand{\byauthors}[1]{#1 }
\newcommand{\journal}[1]{ #1 }
\newcommand{\reftitle}[1]{{\it #1} }
\newcommand{\volumin}[1]{{\bf #1} }
\newcommand{\eref}




\bibitem{1}\byauthors{Opanowicz 1998}
\journal{$J. Phys. A: Math. Gen.$}\volumin{31}{2887.}\eref
\bibitem{2}\byauthors{Kowalczyk 1998 q-alg/9804055}
\journal{}\volumin{}{} \eref
\bibitem{3}\byauthors{Ma\'slanka 1994}
\journal{$J. Math. Phys. A: Math. Gen.$}\volumin{35}{1976.} \eref
\end{thebibliography}
\end{document}